\documentclass[12pt]{amsart}
\usepackage{amsmath}
\usepackage{amssymb}
\usepackage{latexsym}
\usepackage{amsfonts}
\usepackage{amscd}
\usepackage[all]{xypic}
\theoremstyle{plain}
\newtheorem{theorem}{Theorem}[section]

\newtheorem{corollary}[theorem]{Corollary}

\newtheorem{example}[theorem]{Example}

\newtheorem{lemma}[theorem]{Lemma}

\newtheorem{proposition}[theorem]{Proposition}
\newtheorem{remark}[theorem]{Remark}

\newcommand{\N}{{\mathbb N}}

\newcommand{\Z}{{\mathbb Z}}

\newcommand{\gp}{{\pi}}
\newcommand{\gep}{{\epsilon}}
\newcommand{\ga}{{\alpha}}
\newcommand{\gb}{{\beta}}
\newcommand{\irrng}{\textrm{Irr}_n(G)}

\newcommand{\irrngxx}{\textrm{Irr}_n(G^{xx})}

\newcommand{\map}{\textrm{MAP}\ }
\newcommand{\tak}{\textrm{TAK}\ }

\newcommand{\mR}{\mathcal{R}}

\newcommand{\lra}{\longrightarrow}

\newcommand{\nbd}{neighborhood \ }
\begin{document}
\title{Some new results on the Chu duality of discrete groups}
\author{Salvador Hern\'andez and Ta-Sun Wu}
\address{\noindent Departamento de Matem\'aticas,\newline \indent
\'Area Cient\'{\i}fico-T\'ecnica, \newline \indent Universidad
Jaume I,\newline \indent 8029~-~AP Castell\'on, \newline \indent
Spain.} \email{hernande@mat.uji.es}
\address{\noindent Department of Mathematics, \newline
\indent Case Western Reserve University, \newline \indent
Cleveland, Ohio 44106, \newline \indent U.S.A.}
\email{txw3@po.cwru.edu}
\thanks{The first named author acknowledges partial
financial support by the Spanish Ministry of Science (including
FEDER funds), grant MTM2004-07665-C02-01; and the Generalitat
Valenciana, grant GV04B-019 }
\date{}
 \dedicatory{}
\begin{abstract}
This paper deals mainly with the Chu duality of discrete groups.
Among other results, we give sufficient conditions for an $FC$
group to satisfy Chu duality and characterize when the Chu
quasi-dual and the Takahashi quasi-dual of a group $G$ coincide.
As a consequence, it follows that when $G$ is a weak sum of a
family of finite simple groups, if the exponent of the groups in
the family is bounded then $G$ satisfies Chu duality; on the other
hand, if the exponent of the groups goes to infinite then the Chu
quasi-dual of $G$ coincide with its Takahashi quasi-dual. We also
present examples of discrete groups whose Chu quasi-duals are not
locally compact and examples of discrete Chu reflexive groups
which contain non-trivial sequences converging in the Bohr
topology of the groups. Our results systematize some previous work
and answer some open questions in the subject
\cite{chu,re_tri:moore,comfhernremustrig}.
\end{abstract}
\maketitle
\section{Introduction}
It is a consequence of the celebrated Gel'fand and Ra\v{\i}kov
Theorem that the set of all unitary representations of a locally
compact group $G$ contains the information necessary to recover
the topological and algebraic structure of the group (see
\cite{enoschw}). However, in general, such representations may not
be finite dimensional. The result is that if we form a {\it dual
space} associated to the set of all representations, we obtain an
object whose structure is very involved. In \cite{chu} Chu,
motivated by the work of Hochschild and Mostow \cite{hos_most} for
compact groups, considered the groups with enough finite
dimensional representations to separate the points (the so-called
maximally almost periodic groups) and established a duality theory
within this class that extends both Pontryagin and
Tannaka-Kre\v{\i}n dualities for locally compact Abelian groups
and compact groups respectively. Chu defined the space $G^x$
consisting of all finite dimensional representations on $G$
equipped with the compact open topology. Then he used the dual
structure of $G^x$ inherited from $G$ to form a bidual object
$G^{xx}$ consisting of certain continuous mappings on $G^x$ that
behave nicely with respect to the algebraic structure of $G^x$. It
turns out that $G^{xx}$ is a topological group with nice features,
which is called the {\it Chu quasi-dual group} of $G$. When $G$ is
topologically isomorphic to $G^{xx}$ it is said that $G$ satisfies
{\it Chu (unitary) duality} or $G$ is {\it Chu reflexive} (is {\it
Chu} for short). Thus, investigating the Chu duality of locally
compact groups is equivalent to identifying the locally compact
groups that can be recovered from their set of finite dimensional
unitary representations. Chu set the basis for this duality theory
but, after Chu's seminal paper, there have been several important
contributions devoted to developing this theory. We shall mention
some of them below. Nevertheless, many important questions along
this line of research are still open. Some of them were asked by
Chu himself. We mention just a representative one here. Does the
discrete free group with two generators satisfy Chu duality?

This paper deals mainly with the Chu duality of discrete groups.
Among other results, we give sufficient conditions for an $FC$
group to satisfy Chu duality and characterize when the Chu
quasi-dual and the Takahashi quasi-dual of a group $G$ coincide.
As a consequence, it follows that when $G$ is a weak sum of a
family of finite simple groups, if the exponent of the groups in
the family is bounded then $G$ satisfies Chu duality; on the other
hand, if the exponent of the groups goes to infinite then the Chu
quasi-dual of $G$ coincide with its Takahashi quasi-dual. We also
present examples of discrete groups whose Chu quasi-duals are not
locally compact and examples of discrete Chu reflexive groups
which contain non-trivial sequences converging in the Bohr
topology of the groups. Our results systematize some previous work
and answer some open questions in the subject (see
\cite{chu,re_tri:moore,comfhernremustrig}).

\section{Basic Definitions and Facts}
In principle, all groups are assumed to be locally compact
Hausdorff and maxi\-mally almost periodic (MAP). That is, locally
compact groups that can be continuously injected into compact
groups. 
For any group $G$, the symbol $G^{\prime }$ means its commutator
subgroup. 
The group $G$ has {\em finite exponent} if the orders of the
elements of $G$ are finite and bounded. The {\em exponent} of $G$,
denoted $exp(G)$ is the least common multiple of all the orders of
elements in $G$. 

Let $G$ be a topological group, denote by $G_n^x$ the set of all
continuous $\mbox{$n$-dimensional}$ unitary representations of
$G$, i.~e., the set of all continuous homomorphisms of $G$ into
the unitary group $U(n)$, equipped with the compact-open topology.
It follows from a result due to Goto \cite{goto} that $G_n^x$ is a
locally compact and uniformizable space. The space $G^{x}=\sqcup
_{n<\omega }G_n^x$ (as a topological sum) is called the
\textit{Chu dual} of $G$ \cite{chu}.

We now recall the basic notions of the Chu (or unitary) duality.
Its main feature is the construction of a \emph{bidual} of $G$
from the Chu dual $G^x$. This bidual consists of the so-called
quasi-representations. If we define $\mathcal{U}=\sqcup _{n<\omega
}\mathcal{U}(n)$ (topological sum), a
\textit{quasi-representation} of $G^x$ is a continuous mapping
$p:G^x\longrightarrow \mathcal{U}$ conserving the main operations
between unitary representations: direct sums, tensor products,
unitary equivalence and sending the elements of $G_n^x$ into
$U(n)$ for all $n\in \N$ (see \cite{chu} or \cite{heyer70} for
details). The set of all quasi-representations of $G$ equipped
with the compact-open topology is a topological group with
pointwise multiplication as the composition law, called the
\textit{Chu quasi-dual group} of $G$ and denoted by $G^{xx}$.
Thus, a \nbd base of the identity in $G^{xx}$ consists of sets of
the form $[K_n,V]=\{p\in G^{xx} : p(K_n)\subset V \}$, where $V$
is any \nbd of the identity in $U(n)$ and $K_n$ is any compact
subset of $G_n^x$, $n\in \N$. It is easily verified that the
evaluation map $\epsilon_{G} :G\longrightarrow G^{xx}$ is a group
homomorphism which is a monomorphism if and only if $G$ is MAP.
The group $G$ \textit{satisfies Chu duality}
when the evaluation map $\epsilon_{G} $ is an isomorphism of
topological groups. If the evaluation map is only an algebraic
isomorphism, we say that $G$ is \textit{Chu semi-reflexive}.
In this terminology it was shown in \cite{chu}, that LCA groups
and compact groups satisfy Chu duality (indeed Chu duality reduces
to Pontryagin duality and to Tannaka duality respectively for such
groups). Here, the group $\epsilon_G(G)$ is always assumed to be
equipped with the
topology inherited from $G^{xx}$. 

Two $n$-dimensional representations $D_1$ and $D_2$ of $G$ are
said (unitarily) {\it equivalent} to each other when there is $U$
in $U(n)$ such that $D_1(x)=U^{-1}D_2(x)U$ for all $x\in G$. This
notion sets an equivalence relation on $G^x_n$ that we denote with
the symbol "$\sim$".

Let $\irrng$ denote the set of all irreducible representations of
dimension $n$. In addition to the unitary dual, it will be useful
to consider the set $\widehat{G}_n=\irrng/~\sim$ , \ which is the
quotient space formed by the equivalence classes of irreducible
representations of dimension $n$.

The Bohr compactification of an  arbitrary topological   group can
be defined as a pair $(bG,b)$ where  $bG$ is a compact Hausdorff
group  and $b$ is a continuous  homomorphism
 from $G$ onto a dense subgroup of $bG$
with the following universal property: for every continuous
homomorphism $h$ from $G$ into a compact group $K$ there is a
continuous homomorphism $h^{b}$ from $bG$ into $K$ extending $h$
in the sense that $h=h^b \circ b$, that is,
 making the following diagram commutative:
\[
\xymatrix{ G \ar@{>}[rr]^b \ar[dr]^h & & bG\ar[dl]_{h^b} \\ & K &}
\]
  The group $bG$ is essentially
unique and  is also referred to as the   Bohr compactification of
$G$. Heyer \cite[V, \S 14]{heyer70} contains a careful examination
of $bG$ and its properties. The group $G$ may be equipped with the
topology induced by the above homomorphism $b$, the so-called
\textit{Bohr topology}. The Bohr topology is Hausdorff precisely
when $G$ is maximally almost periodic (\textit{MAP} group),
equivalently, when $b$ is one-to-one. Here, we will be mainly
concerned with this class of groups; these turn out to be the
groups whose finite dimensional representations separate points.
The Bohr topology of a group $G$ may also be defined as the one
that $G$ inherits from $\mathcal U^{G^x}$ (the topology of
pointwise convergence on $G^x$). The symbol $G^{b}$ stands for $G$
equipped with the Bohr topology. Clearly, each member of $G^x$
defines a continuous mapping on $G^{b}$ that extends to $bG$.
Thus, the representation spaces $G^x$ and $(bG)^x$ have exactly
the same underlying set. As a consequence, the Chu quasi-dual
group $G^{xx}$ is always algebraically embedded into $bG$.
\medskip

 We now collect some well-known facts about Chu
duality theory (see \cite{heyer70,heyer73,riggins,roeder})

\begin{proposition} \label{pr11}
The following assertions hold for any locally compact MAP group
$G$.
\begin{enumerate}
\item The space $G_n^x$ is locally compact for each $n\in \N$, so
is $G^x$. If $G$ is discrete (resp. metrizable), then $G_n^x$ is
compact (resp. hemicompact) and therefore $G^x$ is
$\sigma$-compact.

\item If $G$ is compact then the quotient space $G^x/\sim$ is
discrete.

\item If $G$ satisfies the second axiom of countability then
$G_n^x$ satisfies the second axiom of countability. So $G^x$ is
metrizable.

\item The group $G$ is maximally almost periodic if and only if
$\epsilon_{G}$ is injective.

\item $G^{xx}$ is complete with respect to uniformity of uniform
convergence on compact subsets.

\item If $G$ is second countable then $G^{xx}$ is a second
countable complete metric space. As a consequence, if the
evaluation map $\epsilon_{G}$ is onto then $G$ is Chu.
\end{enumerate}
\end{proposition}

In some cases, the Chu quasi-dual of a group $G$ coincides with
its Bohr compactification $bG$. 
Here on, we denote by $C(K,U(n))$ the space of all continuous
functions from a topological space $K$ into $U(n)$.
\begin{proposition}\label{discrete}
Let $G$  be a locally compact group. It holds that $G^{xx}$ is
topologically isomorphic to $bG$ if and only if $\widehat{G}_n$ is
discrete for any $n\in \N$.
\end{proposition}
\begin{proof}
{\it Sufficiency:} This is Theorem 4.4 of \cite{gal_her:i0}.

{\it Necessity:} Let $\alpha : G_n^x\longrightarrow
(G^{xx})_{n}^x$ the canonical evaluation map defined by
$\alpha(D)(p)=p(D)$ for all $D\in G_n^x$ and $p\in G^{xx}$. We
first check that $\alpha$ is continuous and injective. For the
continuity, let $\{D_j \}$ be a net converging to $D$ in $G^x_n$.
If $K$ is a compact subset of $G^{xx}$, since $G^x_n$ is locally
compact, it follows that $K$ is an equicontinuous subset of
$C(G^x_n,U(n))$. Hence, the net $\{\alpha(D_j) \}$ converges to
$\alpha(D)$ uniformly on $K$. This yields the continuity of
$\alpha$. For the injectivity, suppose that $D,E$ belong to
$G^x_n$ and $\alpha(D)=\alpha(E)$. Then $p(D)=p(E)$ for all $p\in
G^{xx}$. In particular, $D(g)=E(g)$ for all $g\in G$. Thus, D=E.

Now, we define $\overline{\alpha} : G^x_n\diagup_{ \sim} \
\longrightarrow (G^{xx})_n^x \diagup_{\sim}$ \ so that the
following diagram commutes
\[
\begin{CD}
G_n^x@>\alpha>>   (G^{xx})_n^x\\
@VQVV   @VVQ^{xxx}V\\     
G^x_n\diagup_{\sim} @>\overline{\alpha}>>  (G^{xx})_n^x\diagup_{\sim}
\end{CD}
\]
\medskip

\noindent where $Q$ and $Q^{xxx}$ are the canonical quotient
mappings. In order to check that $\overline{\alpha}$ is properly
defined, suppose that $E$ and $D$ are in $G^x_n$ and $E\sim D$.
Then, there is $A\in U(n)$ such that $D(g)=AE(g)A^{-1}$ for all
$g\in G$. Since, by hypothesis, $G^{xx}$ is topologically
isomorphic to $bG$, it follows that $G$ is dense in $G^{xx}$ with
respect to the topology of pointwise convergence on $G^x$. Thus,
$D(p)=AE(p)A^{-1}$ for all $p\in G^{xx}$. Therefore,
$\alpha(E)\sim \alpha(D)$ and this implies that
$\overline{\alpha}$ is well defined. The fact that $Q$ and
$Q^{xxx}$ are quotient mappings yields the continuity of
$\overline{\alpha}$. Finally, let us see that $\overline{\alpha}$
is also injective. Suppose $E$ and $D$ in $G^x_n$ such that
$\overline{\alpha}(Q(E))=\overline{\alpha}(Q(D))$. Then there is
$A\in U(n)$ such that $\alpha(D)(p)=A\alpha(E)(p)A^{-1}$ for all
$p\in G^{xx}$. Thus, $D(g)=AE(g)A^{-1}$  for all $g\in G$ and, as
a consequence, $E\sim D$. This gives the injectivity of
$\overline{\alpha}$. On the other hand, since $G^{xx}$ is compact,
it follows from Proposition \ref{pr11}(2) that
$(G^{xx})^x_n\diagup_{\sim}$ is discrete. Since
$\overline{\alpha}$ is continuous and injective, it follows that
$G^x_n\diagup_{\sim}$ is also discrete. This completes the proof
because $\widehat{G}_n$ is a subspace of the latter.
\end{proof}
\medskip

We notice that in general the surjectivity of $\epsilon_{G} :
G\longrightarrow G^{xx}$ does not imply that $\epsilon_{G}$ is a
homeomorphism. As an example (due to Moran) let $(p_n)$ be an
infinite sequence of distinct primes numbers ($p_n>2$), and let
$G_n$ be the projective special linear group of dimension two over
the Galois field $(GF(p_n))$ of order $p_n$. Let
$G=\prod_{n=1}^{\infty}G_n$ with discrete topology. Then
$G^{xx}=\prod_{n=1}^{\infty}G_n$ is compact (with product
topology), (cf. \cite{moran}).
\section{$FC$ Groups}

Let $G$ be a (discrete) group. We say that $G$ is an \emph{FC}
group if every conjugacy class of $G$ is finite, i.e., for all
$g\in G$, we have that $\mathcal{O}_x:= \{hgh^{-1}: h\in G \}$ is
finite. We say that $g\in G$ is \emph{central} in $G$ if $gh=hg$
for all $h\in G$. The set of all central elements of $G$ is a
normal subgroup called the \emph{center}, $Z(G)$, of $G$. A
subgroup $H$ of $G$ is called central in $G$ if $H\subset Z(G)$.
In the sequel, the symbol $mdus(G)$ denotes the minimal natural
number $n$ such that the unitary representations of dimension less
or equal than $n$ separate the points in $G$. Our main goal in
this section is to prove the following result.

\begin{theorem}
\label{th_fc1} Let $G$ be an $FC$ torsion group and suppose there
is $N\in \N$ such that $exp(G')\le N$ and $mdus(G/H)\le N$ for all
normal subgroup $H$ of $G$ that is co-finite in $G'$. Then the
group $G$ is Chu reflexive.
\end{theorem}

As a corollary, it follows the following result about direct sums
which are not necessarily torsion.

\begin{corollary}
\label{th21} 
Let $G=\sum_{i\in I}F_{i}$ equipped with the discrete topology
such that $F'_i$ is finite for all $i\in I$. Suppose further that
there is $N\in \N$ with $mdus(F_i)\leq N$ and
$exp(F_i\hspace{0.5mm} ')\leq N$, for all $i\in I$. Then the group
$G$ is Chu reflexive.
\end{corollary}

We notice that the constraints on $mdus(G)$ and $exp(G')$ are only
needed in Proposition \ref{fc_chucommutator} below. The proof of
Theorem \ref{th_fc1} is split in several partial results. Firstly,
we need some definitions.
\medskip

Let $G$ be a topological group and let $H$ be a subgroup of it. We
define the {\it Chu quasi-dual of $H$ in $G$} to be the group
$G^{xx}\cap cl_{bG}H$ 
We denote this group by $(H,G)^{xx}$.

We say that $H$ is {\it Chu semi-reflexive in} $G$ when $H$
coincides
algebraically with $(H,G)^{xx}$.


\begin{proposition}\label{chusubgroup} Let $G$ be a topological
group  and let $H$ be\ 
normal subgroup which is Chu semi-reflexive in $G$. If $G/H$ is
Chu semi-reflexive then $G$ is Chu semi-reflexive.
\end{proposition}
\begin{proof} Consider the exact sequence

\[
1\longrightarrow H\overset{i}{\longrightarrow }G\overset{j}{%
\longrightarrow }G/H\longrightarrow 1
\]

\noindent and the dual maps (in the category of pointed spaces)

\[
1\longrightarrow (G/H)^{x}\overset{j^{x}}{\longrightarrow }
G^{x}\overset{i^{x}}{\longrightarrow }G^{x}_{|i(H)}
\longrightarrow 1
\]
\bigskip

\noindent where $i^{x}$ and $j^{x}$ are both continuous, the map
$i^{x}$ is onto, and the inverse image under $i^{x}$ of the set of
representations of $G$ which are trivial on $H$, is $(G/H)^x$.





Repeating the process above, one can also obtain the Chu
quasi-dual sequence.
\[
1\longrightarrow (H,G)^{xx}\overset{i^{xx}}{\longrightarrow
}G^{xx} \overset{j^{xx}}{\longrightarrow
}(G/H)^{xx}\longrightarrow 1
\]

\noindent Next we verify that this sequence is exact.

Firstly, observe that, since $\epsilon$ is a natural
transformation, the following diagram is commutative

\[
\begin{CD}
1@>>>   H@>i>>             G@>j>>          G/H@>>>            1\\
 @VVV   @V(\epsilon_{G})_{|H}VV     @V\epsilon_{G}VV  @VV\epsilon_{(G/H)}V @VVV\\
1@>>>   (H,G)^{xx}@>>i^{xx}> G^{xx}@>>j^{xx}>(G/H)^{xx}@>>> 1
\end{CD}
\]
\noindent where $(\epsilon_G)_{|H}$ and $\epsilon_{(G/H)}$ are
known to be onto algebraic isomorphisms and $\epsilon_G$ is
1-to-1.

Being and embedding, the map $i^{xx}$ is 1-to-1 and since
$\epsilon_{(G/H)}$ is bijective and the diagram above is
commutative, it follows that that $j^{xx}$ is onto.

Next we show that $Im(i^{xx})=\ker \ j^{xx}$. Indeed,
$i^{xx}(H,G)^{xx}=(\epsilon_G\circ i\circ
\epsilon^{-1}_H)(H,G)^{xx}=\epsilon_G(i(H))\subset \epsilon_G(ker
\ j)\subset ker \ j^{xx}$. On the other hand, let $p$ be an
arbitrary element of $ker \ j^{xx}$. We have that
$j^{xx}(p)(E)=p(j^{x}(E))=p(E\circ j)=I_n$ for all $E\in
(G/H)^{x}_n$ and for all $n<\omega$. We must check that $p$ is in
$cl_{bG}H= K$. Suppose the contrary and let $\pi:bG\longrightarrow
bG/K$ the canonical quotient mapping. We have that $\pi(p)\not=
1_{bG/K}$. Hence, there is $E\in (bG/K)^{x}_n$ such that
$E(\pi(p))\not= I_n$. Now, it is readily seen that $b(G/H)$ is
topologically isomorphic to $bG/K$. Thus, $E_{|(G/H)}$ belongs to
$(G/H)^{x}_n$ and we have that
$p(j^{x}(E_{|(G/H)}))=p(E_{|(G/H)}\circ j)=p(E\circ
\pi)=E(\pi(p))\not= I_n$. This is a contradiction which completes
the proof 
of the exactness of the quasi-dual sequence.

 Applying the well-known five-lemma to the diagram above, using
that $\epsilon_{H}$ and $\epsilon_{(G/H)}$ are isomorphisms onto,
we obtain that $\epsilon_G$ is an algebraic onto isomorphism.
\end{proof}
\begin{corollary}\label{co21} Let $G$ be a LC topological group. If G'\ 
is Chu semi-reflexive in $G$ then the group $G$ is Chu
semi-reflexive.
\end{corollary}

Next lemma is due to Wu and Riggins (see \cite[p.
462]{wu_riggins}) but we include a short proof of it for the
reader's sake. We recall that when $G$ is an FC group we have that
$\mathcal{O}_x$ is finite for all $x\in G$. Thus, $F_x$, the
isotropy group of $x$, is co-finite in $G$.
\begin{lemma}\label{le_fc_central1}
Let $G$ be an $FC$ group. Given an arbitrary element $D$ of
$G^x_n$, $n\in \N$, we have $Z(D(G))$ is co-finite in $D(G)$.
\end{lemma}
\begin{proof}
Let $K$ be the closure of $D(G)$ in U(n). Then $K$ is a compact
Lie group with a dense $FC$ subgroup. Thus, $\mathcal{O}_x$ is
finite for all $x\in D(G)$. In particular, if $K_0$ denotes the
connected component of $K$, we have $\mathcal{O}_x\cap K_0=\{x\}$
for all $x\in D(G)\cap K_0$. Since $D(G)$ is dense in $K$, this
means that $K_0\subset Z(K)$. From the compactness of $K$, it
follows that $K_0$ is co-finite in $K$ and, as a consequence, that
$Z(D(G))$ is co-finite in $D(G)$.
\end{proof}

\begin{lemma}\label{le_fc1}
Let $G$ be an FC group. Given an arbitrary element $D$ of $G^x_n$,
$n\in \N$, there is a co-finite normal subgroup $N$ of $G$ such
that $D(N')=\{I_n\}$. Moreover, the group $G'/N'$ is finite.
\end{lemma}
\begin{proof}
Since $G$ is an FC group, by the lemma above, we have that
$L=Z(D(G))$ is co-finite in $D(G)$. Thus,
$D(G)/L=\{D(x_1)L,...,D(x_m)L\}$. For each $x_i$, $1\leq i\leq m$,
let $F_i$ be the isotropy group of $x_i$. The subgroup
$F=\cap_{i=1}^m F_i$ is co-finite and $ax_ia^{-1}=x_i$ for $a\in
F$, $1\leq i\leq m$. Hence $D(F)\subset Z(D(G))$, what yields
$D(F')=\{I_n\}$. Now, if we take $A$ to be the normal subgroup
generated by $\{x_1,...x_m\}$, then its centralizer in $G$, the
subgroup $C_G(A)$, is a co-finite normal subgroup of $G$ with
$F\supset C_G(A)$ (see \cite[p. 462]{wu_riggins}). Thus, it
suffices to take $N=C_G(A)$ and the first part of proof is done.
Let us verify now that $G'/N'$ is finite. Indeed, since $N$ is
co-finite in $G$, we have $\frac{G}{N}=\{x_1N,...,x_lN\}$, with
$x_i\in G$, $1\leq i\leq l$. On the other hand, $G'$ is generated
by elements of the form $[x_ia,x_jb]$, with $a,b\in N$. Now,
because $N$ is normal in $G$, we know that $N'$ is also normal in
$G$. Hence, we have
\begin{equation*}
\begin{array}{l}
[x_ia,x_jb]N'=\\
(x_iax_jba^{-1}x_1^{-1}b^{-1}x_j^{-1})N'=\\
(x_i(ax_ja^{-1})[a,b](bx_i^{-1}b^{-1})x_j^{-1})N'=\\
(x_i(ax_ja^{-1})(bx_i^{-1}b^{-1})x_j^{-1})N'
\end{array}
\end{equation*}

Since $G$ is an FC group, we obtain that there are only finitely
many elements of the form $[x_ia,x_jb]N'$, with $a,b\in G$ and
$x_i,x_j\in \{x_1,...,x_l\}$, which proves that $\frac{G'}{N'}$ is
finitely generated. Now, by \cite[14.5.9]{robinson}, $G'$ is
torsion and, as a consequence, so is $\frac{G'}{N'}$. Hence, the
quotient group $\frac{G'}{N'}$ is finitely generated, torsion and
FC. Applying \cite[14.5.7]{robinson}, we obtain that
$\frac{G'}{N'}$ is finite, which completes the proof.
\end{proof}

Using that $N'\subset Ker D$ in the lemma above, it follows.
\begin{corollary}\label{co_fc1}
Let $G$ be an FC group. Then for any element $D$ of $G^x_n$, $n\in
\N$, we have that $G'/((Ker D)\cap G')$ is finite.
\end{corollary}

The proof of Lemma \ref{le_fc1} may be adapted easily to obtain
also following result. 
\begin{lemma}\label{riggins}
Let $G = \sum_{i \in I} F_i$, where $F'_i$ is finite for all $i\in
I$. Then, for each $D\in Hom(G, U(n))$, there is a finite subset
$J \subseteq I$ such that, if $x$ and $y$ belong to $G'$ and
$x_{|J} = y_{|J}$, then we have that $D(x) = D(y)$.
\end{lemma}

Let $G$ be an $FC$ group and consider a family $\mathcal B$ of
normal subgroups of $G$ satisfying: (1) for every $A\in \mathcal
B$ we have that $G'/(A\cap G')$ is finite; 
(2) for all finite subset $F$ of $G$, and $D\in G^x_n$ with
$D(F)=\{I_n\}$, there is $A\in \mathcal B$ such that $F\subset A$
and $A\subset Ker D$; and (3) the set $\mathcal B$ is directed
under inverse inclusion. For instance, we have just verified that,
if $\mathcal D\subset G^x$ is such that $\mathcal D$ contains
$G^x_1$ and separates the points of any quotient $G/L$, with $L$ a
normal subgroup of $G$ such that $L\cap G'$ is co-finite in $G'$,
then one can take $\mathcal B$ to be the family $\{\cap_{i=1}^n
Ker D_i: D_i\in \mathcal D, n\in \N \}$. Now, consider the
projective system $\{(G/A,f_{AB}): A,B\in \mathcal B \}$, where
$f_{AB}:G/B\longrightarrow G/A$ is the canonical homomorphism
naturally defined when $B\subset A$ and equip each group $G/A$
with the weak topology generated by $(G/A)^x$. Define $H$ to be
the projective limit $\underleftarrow{Lim} \ (G/A)$, where $f_A$
denotes the canonical homomorphism of $H$ onto $G/A$ and observe
that, since its topology is defined by finite dimensional unitary
representations, the group $H$ is totally bounded. On the other
hand, according to the definition of $\mathcal B$, the subgroup of
$H$, defined as $K=\underleftarrow{Lim} \ (G'/(A\cap G'))$ is
profinite and, therefore, compact.

\begin{proposition}\label{pr_fc1}
Let $G$ a maximally almost periodic $FC$ group and let $\mathcal
B$ be a family of subgroups of $G$ satisfying the assertions (1)-
(3), as defined above. Then $cl_{bG}G'$ is topologically
isomorphic to $K=\underleftarrow{Lim} \ (G'/(A\cap G'))$. 
\end{proposition}
\begin{proof}
Let $\pi_A:G\longrightarrow G/A$ the canonical quotient and
consider the map $\psi: G \rightarrow H=\underleftarrow{Lim} \
(G/A)$, which is defined by $f_A\circ \psi =\pi_A$. Then the map
$\psi$ is a homomorphism trivially continuous when we consider the
Bohr topology on $G$. By the universal property of the Bohr
compactification, $\psi$ can be extended to a continuous group
homomorphism $\bar \psi: bG \rightarrow \overline{H}$, with $\bar
H$ being the Weil completion of the group $H$. Since
$\psi(G^\prime)$ is dense in the compact subgroup
$K=\underleftarrow{Lim} \ (G'/(A\cap G'))$ of $H$, it follows that
$\bar \psi$ is a quotient homomorphism from $cl_{bG} G^\prime$
onto $K$. In order to prove that $\bar\psi$ is an isomorphism, it
suffices to show that $\bar \psi$ is injective on $cl_{bG}
G^\prime$. Now, if $g \in cl_{bG} G^\prime$ and $g \neq 1$, then
there exists $n \in \mathbb{N}$ and $D \in (bG)_n^x=G_n^x$ such
$D(g) \neq I_n$. By Corollary \ref{co_fc1}, we know that the group
$G'/((Ker D)\cap G')$ is finite. By (2), we can take $A\in
\mathcal B$ such that $A\subset G\cap Ker D$. Clearly, the
representation $D$ factorizes through
$E:bG/(cl_{bG}A)\longrightarrow U(n)$. Thus, if $\overline{f}_A$
denotes the canonical extension of $f_A$, the following diagram
\[
\xymatrix{  bG \ar@{>}[rr]^{\overline{\psi}} \ar@{>}[d]_{D} \ar@{>}[drr]^{\overline{\pi}_A} & & \overline{H}
\ar@{>}[d]^{\overline{f}_A} \\
 U(n)  &  & \frac{bG}{(cl_{bG}A)} \ar@{>}[ll]^{E}}
\]
commutes.\medskip

\noindent This yields $\overline{f}_A(\overline{\psi}(g))\not= 1$
and shows the injectivity of $\bar \psi$, 
which completes the proof.
\end{proof}

As a consequence of Proposition \ref{pr_fc1}, we obtain.
\begin{corollary}\label{co25}
Let $G = \sum_{i \in I} F_i$, where $F'_i$ is finite for all $i\in
I$. Then we have that $cl_{bG}G'$ is topologically isomorphic to
$K=\prod_{i\in I}F_{i}^{\prime }$. As a consequence, it follows
that $G'$ is Bohr closed in $G$.
\end{corollary}
\begin{proof}
It is enough to take $\mathcal B = \{\sum_{i\in I\setminus
J}F_{i}: J\subset I,\ \hbox{and}\ J\ \hbox{ finite} \}$
\end{proof}
\begin{proposition}\label{fc_chucommutator}
Let $G$ be an $FC$ torsion group and suppose there exist a family
$\mathcal B$ of normal subgroups of $G$ satisfying (1)-(3), as
defined above. Assume further that there is $N\in \N$ such that
$exp(G')\le N$ and 
$mdus(G/B)\le N$ for all $B\in \mathcal B$. Then $G'$ is Chu
semi-reflexive in $G$.
\end{proposition}
\begin{proof}
It suffices to verify that $G'$ is Chu semi-reflexive in $G$. Let
$g=(g_A)$ be an element of $cl_{bG}G'\cap G^{xx}\cong
\underleftarrow{Lim} \ (G'/(A\cap G'))\cap G^{xx}$ and suppose
that $g\notin G'$. Firstly, observe that $g\notin G$ since $G'$ is
Bohr closed in $G$ (indeed, if $g\in G\setminus G'$, then there is
a character $\chi$ on $G$ such that $\chi(g)\not= 1$). Now, for
any
arbitrary finite subset $F$ of $G$,\ 
since $G$ is an $FC$ torsion group, we know that the normal
closure $F^G$ of $F$ in $G$ is finite. Thus, the same happens with
its normal closure in $bG$. That is, we have $F^{bG}=F^G$. Now,
the compact group $bG/F^G$ is MAP and, since $g\notin F^G$, we can
find $D\in (bG)^x_n$ such that $D(F^G)=\{I_n\}$ and $D(g)\not=
I_n$. By (2), take $B_F\in \mathcal B$ such that $F^G\subset B_F$
and $B_F\subset (G\cap Ker D)$. It follows that
$f_{B_F}(\psi(F^G))=1$ but $\overline{f}_{B_F}(\bar \psi(g))\not=
1$. Therefore, by hypothesis, there is $\tilde{E}_F\in
(G/B_F)^x_N$ such that $\tilde{E}_F(\overline{f}_{B_F}(\bar
\psi(g))\not=I_N$.\
Let $V$ be a \nbd of the identity in $U(N)$ that contains no
non-trivial elements of order $\leq N$. Since $exp(G'/B_F)\le N$
and $\bar\psi_{B_F}(g)\in G'/B_F$,  we have $\tilde{E}_F(\bar
\psi_{B_F}(g)\notin V$.
Set $E_F=\tilde{E}_F\circ \psi_{B_F}$ for all finite subset $F$ of
$G$ and let $\mathcal F$ be the set consisting of all finite
subsets of $G$, ordered under inclusion. Since $E_F(F^G)=\{I_N\}$
for all $F\in \mathcal F$, it is easily verified that the net
$\{E_F: F\in \mathcal F \}$ converges to the identity
representation in $G^x_N$. On the other hand, since $g\in
(\underleftarrow{Lim} \ (G'/A))\cap G^{xx}$, it follows that the
net $\{g(E_F): F\in \mathcal F \}$ must converge to the identity
$I_N$. Nevertheless, $g(E_F)=\tilde{E}_F(\bar\psi_{B_F}(g))\notin
V$ for all $F\in \mathcal F$. This contradiction completes the
proof.
\end{proof}

When $G=\sum_{i\in I} F_i$ is a direct sum, then it is not longer
necessary to assume $G$ to be torsion. Indeed, a slight variation
in the proof of Proposition \ref{fc_chucommutator} also gives the
following result.
\begin{corollary}\label{chucommutator}
Let $G = \sum_{i \in I} F_i$, where $F'_i$ is finite for all $i\in
I$. if we further assume that \ $mdus(F_i)\leq N$ and
$exp(F_i\hspace{0.5mm} ')\leq N$, for some $N\in \N$ and for all
$i\in I$, then $G'$ is Chu semi-reflexive in $G$.
\end{corollary}
\begin{proof}
Again, it is enough to take $\mathcal B = \{\sum_{i\in I\setminus
J}F_{i}: J\subset I,\ \hbox{and}\ J\ \hbox{ finite} \}$
\end{proof}
\begin{lemma}\label{dense}
Let $G$ be a metrizable group such that $\epsilon_G(G)$ is dense
in $G^{xx}$ and the groups $G$ and $\epsilon_G(G)$ have the same
convergent sequences. Then $G^{xx}$ is Chu reflexive.
\end{lemma}
\begin{proof}
It follows from the hypothesis that the group $G$ is MAP. Indeed,
take $g_1 \not= g_2$ in $G$ such that
$\epsilon_G(g_1)=\epsilon_G(g_2)$, then the sequence $\{ x_n \}$
which alternates $g_1$ and $g_2$ is not convergent, but
$\epsilon_G(x_n)$ is convergent which is impossible. 
Since the group $G$ is metrizable, by Proposition 2.1 (1), $G_n^x$
is hemicompact for all $n\in \mathbb{N}$, and thus the
compact-open topology on $G^{xx}$ is metrizable. In particular,
$\epsilon_G(G)$ is metrizable, and thus $G$ and $\epsilon_G(G)$
are topologically isomorphic, because they have the same
convergent sequences. Therefore,  $G^x = (\epsilon_G(G))^x$. On
the other hand, $\epsilon_G(G)$ is dense in $G^{xx}$, so by a
recent result of Luk\'acs (see \cite{lukas}), $(\epsilon_G(G))^x_n
= G^{xxx}_n$. Hence, $G^x = (\epsilon_G(G))^x = G^{xxx}$, and as a
consequence $G^{xx}$ is Chu reflexive.

\end{proof}
\begin{corollary}\label{semi-reflexive}
Let $G$ be a semi-reflexive metrizable group such that $G$ and
$G^{xx}$ have the same convergent sequences. Then $G$ is Chu.
\end{corollary}
\medskip

\begin{proof}[\bf{Proof of Theorem \ref{th_fc1}}]

\noindent \ \ \ Propositions \ref{chusubgroup} and
\ref{fc_chucommutator} yield the semi-reflexivity of $G$. In order
to apply Corollary \ref{semi-reflexive}, we have to verify that
$G$ and $G^{xx}$ have the same convergent sequences. Since the
evaluation map is always continuous, it follows that every
sequence converging in $G$ also converges in $G^{xx}$. So, let
$\{\overline{g}_m \}$ be a sequence converging to the neutral
element in $G^{xx}$. By hypothesis, there must be a co-countable
normal subgroup $L$ of $G$ such that the canonical quotient
$\pi_L:G\longrightarrow G/L$ is one-to-one on the sequence
$\{\overline{g}_m \}$. On the other hand, since every
representation on $G/L$ can be extended canonically to $G$, we
deduce that $\{\pi_L(\overline{g}_m) \}$ converges to the neutral
element in $(G/L)^{xx}$. Now, the group $G/L$ is countable and, by
the same logic as $G$, semi-reflexive . Hence, by Proposition
\ref{pr11} (6), $G/L$ is Chu. That is to say,
$\{\pi_L(\overline{g}_m) \}$ converges in $G/L$ which is a
discrete group. Hence, the injectivity of $\pi_L$ on the sequence
$\{\overline{g}_m \}$ completes the proof.
\end{proof}
\medskip

Theorem \ref{th_fc1} yields the following characterization of Chu
reflexivity for direct sums of finite simple non-abelian groups.

\begin{corollary}\label{co39}
Let $G=\sum_{i\in I}F_i$, where $F_i$ is a finite simple
non-abelian group for each $i\in I$. Then the group $G$ is Chu if
and only if the set $\{exp(F_i) : i\in I \}$ is bounded.
\end{corollary}
\begin{proof}
Assume that the set $\{exp(F_i) : i\in I \}$ is bounded. Let $X$
be a set of finite cardinality, $|X|$, greater than $\exp(F_i)$
for all $i\in I$. By Cayley's Theorem (see
\cite[1.6.8]{robinson}), every group $F_i$ is isomorphic to a
subgroup of $S_{|X|}$, the symmetric group of degree $|X|$. Thus,
the set $\{mdus(F_i) \}$ is bounded. Finally, Corollary \ref{th21}
yields that $G$ is Chu.


Conversely, suppose that $\{exp(F_i) : i\in I \}$ contains a
countable subset, say $\{F_{i_j} : j\in \N \}$, such that
$exp(F_{i_j})\geq j$ for all $j\in \N$. Since the groups $F_{i_j}$
are simple and non-abelian, it follows that if $[G^x_n,V]$ is an
arbitrary canonical \nbd of the identity in $G^{xx}$, for $V$ any
\nbd \ of the identity in $U(n)$, it holds that $[G_n^x,V]$ always
contains a subgroup of the form $\sum_{j\geq m} F_{i_j}$, for some
$m\in \N$. Indeed, 
applying Jordan's theorem (see \cite{curtis_reiner}), given any
$n\in \N$, there is a natural number $f(n)$ such that, for all
$D\in G^x_n$, $D(F_{i_j})=I_n$ if $j\geq f(n)$. This clearly
implies that $[G^x_n,V]$ contains the subgroup $\sum_{j\geq f(n)}
F_{i_j}$ . This verifies that $G^{xx}$ is not discrete. Therefore,
the group $G$ is not Chu.
\end{proof}
\begin{remark}
Observe that, if $C>0$ is given and $\{F_i\}_{i\in I}$ is a family
of finite simple groups, such that $exp(F_i)\leq C$ for all $i\in
I$, then there are finitely many $F_i$'s up to isomorphism: this
follows from the fact that every finite simple group is
$2$-generated, and from the positive solution to the restricted
Burnside problem.
\end{remark}
\medskip

Following the terminology of Trigos-Arrieta, we say that a
topological group $G$ {\it respects} a topological property $P$ if
a subset $A$ of $G$ has $P$ as a subspace of $G$ if and only if
$A$ has $P$ as a subspace of $bG$. In \cite{re_tri:moore} (see
also \cite{comfhernremustrig}) Remus and Trigos-Arrieta ask
whether a Chu group $G$ always respects compactness. We notice
that if $G$ is the weak sum of a countable family of finite simple
non-abelian groups whose order is bounded then $G$ is Chu by
Corollary \ref{co39}. On the other hand, Corollary \ref{co25}
implies that $bG=\prod_{i\in \N} F_i$ which is metrizable. This
means that the discrete group $G$ contains non-trivial sequences
that are convergent in $bG$. Thus, there are Chu groups that do
not respect compactness. \medskip

The following examples follow from Theorem \ref{th_fc1}. Here on,
the symbolism $H\rtimes K$ denotes the semi-direct product of the
groups $H$ and $K$ where the latter is a subgroup of the
automorphism group of the former.
\begin{example}(Heyer \cite{heyer73})
Let $\Z_3\rtimes \Z_2=S_3$ the permutation group. Define
$G_i=\Z_3\rtimes \Z_2$ for all $i\in \N$ and take $G= \sum_{i\in
\N}G_i$. Then $G$ is Chu reflexive.
\end{example}

That not every direct sum of finite groups is Chu was shown by
Moran.
\begin{example}(Moran \cite{moran})
Let $\{p_i\}$ be an infinite sequence of distinct prime numbers
($p_i>2$), and let $F_i$ be the projective special linear group of
dimension two over the Galois field $GF(p_i)$ of order $p_i$. It
holds that $G=\sum_{i\in \N}F_i$ is not Chu.
\end{example}
\medskip

%
%

The methods of Theorem \ref{th_fc1} can also be applied to some
non discrete groups with slight modifications. Next follows an
example of this fact. Let $\{G_i:i\in I\}$ be a non void family of
locally compact groups and set $G=\prod_{i\in I} G_i$. Let $H_i$
be an open subgroup of $G_i$ for all $i\in I$. The {\em local
direct product of the $G_i$'s relative to the open subgroups
$H_i$} is defined as the subgroup $G_0$ of $G$ consisting of all
$(x_i)$ for which $x_i\in H_i$ for all but a finite number of
indices $i$, and equipped with the following topology: Let $J$ be
a finite subset of $I$, and let $U_i$ be a \nbd \ of $1_i$ in the
subgroup $H_i$ for each $i\in J$; let the sets $\prod_{i\in J} U_i
\times \prod_{i\in I\setminus J} H_i$ be taken as an open basis at
the neutral element $(1_i)$ in $G_0$. It follows that $G_0$ is a
topological group containing
the subgroup $H=\prod_{i\in I} H_i$ as an open subgroup. 
It is readily seen that $G_0$ is locally compact if the subgroup
$H_i$ is compact for all $i\in I$. We have the following
consequence of Corollary \ref{th21}.
\begin{corollary}\label{localproduct}
Let $G_0$ be the local direct product of the family $\{G_i,H_i\}$\
where $H_i$ is a compact normal subgroup of $G_i$ such that
$G'_i/H_i$ is finite \
for all $i\in I$. Suppose further that $mdus(G_i/H_i)\leq N$
and $exp(G_i\hspace{0.5mm}'/H_i)\leq N$, for some $N\in \N$ and for all $i\in I$. 
Then $G_0$ is Chu semi-reflexive. Moreover, if $I$ is countable,
the group $G_0$ is Chu reflexive.
\end{corollary}
\begin{proof}
Observe that $H= \prod_{i\in I}  H_i$ is a compact-open subgroup
of $G_0$. Hence, $G_0/H$ is discrete and thus topologically
isomorphic to $\sum_{i\in I} G_i/H_i$. The collection $\{G_i/H_i :
i\in I \}$ satisfies the conditions of Theorem 3.1, so its
(weak-)sum is Chu semi-reflexive. Therefore $G_0/H$ is
semi-reflexive. Since $H$ is compact, it is trivially Chu
semi-reflexive in G. Now, it suffices to apply Proposition
\ref{chusubgroup} to obtain that G is semi-reflexive. Finally, in
case $I$ is countable, the local direct product $\{G_i,H_i \}$ is
second countable and locally compact. Hence, the reflexivity
follows immediately from item (6) of Proposition \ref{pr11}.
\end{proof}
\medskip

Corollary \ref{localproduct} does not hold if we remove the
constraint of being $H_i$ a normal subgroup of $G_i$.
\begin{example} Take the permutation group $S_3=\Z_3\rtimes
\Z_2$, let $G_i=\Z_3\rtimes \Z_2$ and $H_i=\Z_2$ for all $i\in
\N$. If we consider the local direct product $G_0$ of the family
$\{G_i,H_i\}_{i\in \N}$, then $G_0$ is not Chu reflexive.
\end{example}
\begin{proof}
Indeed, using that the smallest normal subgroup that contains
$H_i$ is $G_i$, it can be proved that the spaces
$(G_0)^{x}_{n}/\sim$ are discrete for all $n\in \N$. Then
Proposition \ref{discrete} yields that $G_0^{xx}$ is topologically
isomorphic to $\prod_{i\in \N} G_i$.
\end{proof}

\section{Takahashi groups}

For each locally compact group $G$, Takahashi has constructed a
locally compact group $G_T$ called {\em Takahashi quasi dual}
such that $G_T$ is maximally almost periodic, and $G_T'$ is
compact. The category of locally compact groups with these two
properties is denoted by \tak.
If $n>1$ and $D\in Hom_c(G,U(n))$ then the sets
$t_n(D;U)=\{D\otimes \chi : \chi\in U \}$, $U$ any \nbd of the
identity in the group $G^x_1$, form a fundamental system of
neighbourhoods of $D$ for a topology in $Hom_c(G,U(n))$. We denote
by $G^t_n$ the set $Hom_c(G,U(n))$ equipped with this topology and
the symbol $G^t$ denotes the topological sum of the spaces
$G^t_n$, for $n\in \N$. A \textit{unitary mapping} on $G^t$ is a
continuous mapping $p:G^t\longrightarrow \mathcal{U}$ conserving
the main operations between unitary representations: direct sums,
tensor products, unitary equivalence and sending the the elements
of $G_n^t$ into $U(n)$ for all $n\in \N$ (see \cite{pogunt:tak}
for details). The set of all unitary mappings on $G^t$ equipped
with the compact-open topology is a topological group, with
pointwise multiplication as the composition law, which called the
\textit{Takahashi quasi-dual group} of $G$ and is denoted by
$G_T$. It is easily verified that $G^{xx}\subset G_T\subset bG$.
On the other hand, the evaluation map defines a group homomorphism
$\epsilon^T_{G} :G\longrightarrow G_T$, which is a monomorphism if
and only if $G$ is MAP. The group $G$ \textit{satisfies Takahashi
duality}
when $\epsilon^T_{G} $ is an isomorphism of topological groups.
The Takahashi duality theorem establishes that $G$ satisfies this
duality if $G\in \tak$. It is known that LCA groups and compact
groups belong to \tak. (cf. \cite{taka}). On the other hand, given
any group $G$, we have that $\epsilon_G(G)$ is dense in $G_T$ and
for each $H\in \tak$ and each homomorphism $f:G\lra H$ there
exists exactly one homomorphism $f^t:G_T:\lra H$ with $f=f^t\circ
\epsilon_G$. A detailed discussion and extension of this theory
has been given by Poguntke in \cite{pogunt:tak}, from where we
have taken these lines.

In this section, we are interested in finding out when the Chu
quasi-dual coincides with the Takahashi quasi-dual for \map
groups. Firstly, we give some examples that illustrate the
different situations arising in the theory.

\begin{example}(Moran \cite{moran})
Let $\{p_i\}$ be an infinite sequence of distinct prime numbers
($p_i>2$), and let $F_i$ be the projective special linear group of
dimension two over the Galois field $GF(p_i)$ of order $p_i$. If
$G=\sum_{i\in \N} F_i$, we have $G^{xx}=G_T=bG$.
\end{example}

\begin{example}(Heyer \cite{heyer73})
Let $\Z_3\rtimes \Z_2=S_3$ the permutation group. Define
$G_i=\Z_3\rtimes \Z_2$ for all $i\in \N$ and take $G= \sum_{i\in
\N}G_i$. Then $G=G^{xx}$ and $G_T=\prod_{i\in \N} \Z_3\rtimes
\sum_{i\in \N} \Z_2$.
\begin{proof}
It was proved by Heyer \cite{heyer73} that $G=G^{xx}$. On the
other hand, by Corollary \ref{co25} $cl_{bG}G'=\prod_{i\in \N}
\Z_3$. It suffices now to use \cite[Prop. 17]{pogunt:tak} to
complete the proof.
\end{proof}
\end{example}

\begin{example}
Let $p$ a prime number greater than $2$, and let $F_i$ be the
projective special linear group of dimension two over the Galois
field $GF(p)$ of order $p$ . If $G=\sum_{i\in \N} F_i$, we have
$G^{xx}=G$ and $G_T=bG$.
\begin{proof}
Applying Corollary \ref{th21}, we know that $G^{xx}=G$ and, by
Corollary \ref{co25}, we have $bG=cl_{bG}G'=\prod_{i\in \N} F_i$.
Again, it suffices to apply \cite[Prop. 17]{pogunt:tak} and the
proof is done.
\end{proof}
\end{example}

If $G$ is a locally compact group and $G'$ denotes its commutator
subgroup, we
define the following equivalence relation $\mR$ on 
$\widehat{G}_n$. Let $\sigma, \tau \in \widehat{G}_n$, then
$\sigma \mR \tau$ if and only if there are two representations
$D\in \sigma$, $E\in \tau$ and  an unitary matrix $U\in U(n)$ such
that $D(x)=UE(x)U^{-1}$ for all $x\in G'$. Since each equivalence
class defined by $\sim$ is clearly contained in some equivalence
class defined by $\mR$, it is easily verified that $\mR$ defines a
closed equivalence relation on $\widehat{G}_n$. Hence, there is a
canonical quotient mapping $\overline{\gp}:\widehat{G}_n\lra
\widehat{G}_n/\mR$, which is continuous. On the other hand, we
note that, using essentially the same definition, the equivalence
relation $\mR$ can also be defined on the spaces $\irrng$ and we
obtain a similar quotient mapping $\tilde{\gp}:\irrng\lra
\irrng/\mR$. If we define $p:\frac{\irrng}{\mR}\lra
\frac{\widehat{G}_n}{\mR}$ so that the following diagram
\[
\xymatrix{ \irrng \ar@{>}[rr]^{\gp} \ar@{>}[d]^{\tilde{\gp}} & & \widehat{G}_n \ar@{>}[d]^{\overline{\gp}}\\
\frac{\irrng}{\mR} \ar@{>}[rr]^p & & \frac{\widehat{G}_n}{\mR}}
\]
commutes, then it is easy to verify that the map $p$ is actually a
homeomorphism. Using this observation, we can state now the main
result of this section follows.
\begin{theorem}\label{th_taka_1}
Let $G$  be a locally compact group. We have that $G^{xx}$ is
topologically isomorphic to $G_T$ if and only if
$\widehat{G}_n/\mR$ is discrete for all $n\in \N$.
\end{theorem}

Firstly, we need some preliminary results.

\begin{lemma}\label{le_taka_1}
Let $G$ be a Takahashi group. Then $\widehat{G}_n/\mR$ is discrete
for all $n\in \N$
\end{lemma}
\begin{proof}
According to the definition of $\mR$, we have a commutative
diagram
\[
\xymatrix{ \widehat{G}_n \ar@{>}[rr]^\gp \ar[dr]^r & & \widehat{G}_n/\mR \ar[dl]_{\overline{r}} \\
& (G'_n)^x &}
\]
where, for each $\sigma\in \widehat{G}_n$ and $D\in \sigma$, the
class $r(\sigma)$ is defined by the representation $D_{|G'}\in
Hom(G',U(n))$. On the other hand, the map $\overline{r}$ is
defined by the equality $\overline{r}\circ \gp=r$. It is readily
seen that $\overline{r}$ is continuous and injective. Now, since
$G'$ is a compact group, we have that $(G'_n)^x$ is discrete.
Hence, its inverse image $\widehat{G}_n/\mR$ is also discrete.
\end{proof}

\begin{lemma}\label{le_taka_01}
Let $G$ be a maximally almost periodic, locally compact group such
that $(G^{xx})'$ is a compact subgroup of $G^{xx}$. Then $G/G'$ is
topologically isomorphic to $G^{xx}/(G^{xx})'$.
\end{lemma}
\begin{proof}
Consider the following diagram, where $\gep_G$ is the natural
injection of $G$ into $G^{xx}$ and the other mappings are defined
in canonical way.
\[
\xymatrix{ G^{xx} \ar@{>}[rr]^{\pi^{xx}} \ar[dr]^p & & (G/G')^{xx} \\
&  \frac{G^{xx}}{(G^{xx})'} \ar[ur]_{\beta}  &\\%
 G \ar@{>}[rr]_{\pi} \ar[uu]^{\gep_G} &  &  \frac{G}{G'} \ar[ul]_{\alpha} \ar[uu]_{\theta}}
\]
Let $H$ be the kernel of $\gb$ and let
$\gamma:G^{xx}\longrightarrow bG$ be the natural inclusion (here,
we view $G^{xx}$ as a subgroup of $bG$). We wish to show
$\gamma(H)\subset b(G)'$. First, we note that, since $(bG)'$ is
the intersection of kernels of all unitary characters of $bG$ we
have $\overline{b(G')}=(bG)'$ for every discrete group $G$. Now,
since $(G^{xx})'\supseteq \gep_G(G')$, we have
$\overline{(G^{xx})'}\supseteq \overline{\gep_G(G')}$. Hence,
$\overline{\gamma(\overline{(G^{xx})'})}\supseteq
\overline{\gamma(\gep_G(G'))}=\overline{b(G')}$. As $(G^{xx})'$ is
compact, we obtain $\gamma((G^{xx})')=(bG)'$. Now, for each $h\in
H$, $\beta(h)$ is the identity of $(\frac{G}{G'})^{xx}$; i.e., for
each finite-dimensional unitary representation
$D:\frac{G}{G'}\longrightarrow U(d(D))$, we have
$\gb(h)(D)=I_{d(D)}$. Thus, if $\gp:G\longrightarrow \frac{G}{G'}$
is the canonical quotient mapping, it holds $h(D\circ
\gp)=I_{d(D)}$. Suppose $h\in H$ and $\gamma(h)\notin (bG)'$; then
there is a finite-dimensional unitary representation $E$ of
$bG/(bG)'$ such that $E(\gamma(h)(bG)')\not= I_{d(E)}$. Hence, if
$b_{G'}:\frac{G}{G'}\longrightarrow \frac{bG}{(bG)'}$ denotes the
inclusion canonically associated to $b:G\longrightarrow bG$, we
have $h(E\circ b_{G'}\circ \gp)\not=I_{d(E)}$. This is a
contradiction; therefore, we have $\gamma(H)\subset (bG)'=
\gamma((G^{xx})')$. This implies that the kernel of $\gb$ is
trivial. On the other hand, since $G/G'$ is Abelian, the map
$\theta$ is a surjective topological isomorphism. So, by the
commutativity of the diagram, the map
$\ga:\frac{G}{A}\longrightarrow \frac{G^{xx}}{(G^{xx})'}$ is a
topological isomorphism as well, which completes the proof.
\end{proof}

\begin{proof}[\bf{Proof of Theorem \ref{th_taka_1}}]
Let us suppose that $G_T=G^{xx}$ and consider the diagram
\[
\begin{CD}
\irrng@>\alpha>>   \irrngxx\\
@V\tilde{\gp} VV   @VV\tilde{\gp}^{xx}V\\     
\irrng\diagup \mR @>\overline{\alpha}>> \irrngxx\diagup
\mR
\end{CD}
\]
where $\alpha : \irrng\longrightarrow \irrngxx$ is the canonical
evaluation map defined by $\alpha(D)(x)=x(D)$ for all $D\in
\irrng$ and $x\in G^{xx}$. We know that $\ga$ is continuous and
injective (see Proposition \ref{discrete}). On the other hand, the
map $\overline{\alpha}$ is defined so that the diagram commutes.
Since $G'$ is dense in $G'_T=(G^{xx})'$, it follows that
$\overline{\ga}$ is also well defined, continuous and injective.
Then we set the diagram
\[
\begin{CD}
\irrng\diagup \mR @>\overline{\alpha}>>   \irrngxx\diagup \mR\\
@VpVV   @VVp^{xx}V\\     
\widehat{G}_n\diagup \mR @>\widehat{\alpha}>>
\widehat{G^{xx}}_n\diagup \mR
\end{CD}
\]
where, for any $\sigma\in \widehat{G}_n\diagup \mR$, the element
$\widehat{\ga}(\sigma)$ is defined as the class of
$\overline{\ga}(D)$ with $D\in \sigma$, which makes commutative
the diagram. Again, it is readily verified that $\widehat{\ga}$ is
well defined, continuous and injective. Since $(G^{xx})'$ is
compact, we get by Lemma \ref{le_taka_1} that
$\widehat{G^{xx}}_n\diagup \mR$ is discrete, which implies that
$\widehat{G}_n\diagup \mR$ is discrete.

Conversely, let us suppose that $\widehat{G}_n\diagup \mR$ is
discrete for all $n\in \N$. Since $G$ is dense in $G^{xx}$ with
respect to the Bohr topology (that is to say, when they are
identified to subgroups of $bG$ ), it follows that for every $D$
and $E$ in $Hom_c(G,U(n))$ with $D\mR E$, there is a unitary
matrix $U$ such that $D(x)=UE(x)U^{-1}$ for all $x\in (G^{xx})'$.
Then, repeating an argument similar to the one used in \cite[Th.
4.4]{gal_her:i0}, we obtain that $(G^{xx})'$
is totally bounded as a (topological) subgroup of $G^{xx}$. The
latter group is a complete uniform 
space, which yields $(G^{xx})'$ is compact and topologically
isomorphic to $cl_{bG} G'$. Hence, we have verified that $G^{xx}$
is a Takahashi group. On the other hand, Lemma \ref{le_taka_01}
asserts that $G^{xx}/(G^{xx})'$ is topologically isomorphic to
$G/G'$.\
From this fact, it follows that the canonical mapping
$\epsilon_G:G\lra G^{xx}$ is a dense continuous injection. Now,
consider the diagram
\[
\xymatrix{  G \ar@{>}[rr]^{\epsilon_G} \ar@{=}[d] & & G^{xx} \ar@<1ex>[d]^{id^{xx}} \\
G \ar@{>}[rr]_{\epsilon_G^T} & & G_T \ar@<1ex>[u]^{id^t}}
\]

Since $\epsilon_G$ and $\epsilon^T_G$ are dense continuous
injections into Takahashi groups, we have that $G^{xx}$ and $G_T$
are canonically isomorphic, which completes the proof.
\end{proof}

\begin{corollary}\label{co_taka_1}
Let $G$ be a simple MAP discrete group (which implies $G'=G$).
Then the following conditions are equivalent:
\begin{itemize}
\item[(i)] $G^{xx}=G_T$;

\item[(ii)] $\widehat{G}_n$ is discrete for all $n\in \N$;

\item[(iii)] $G^{xx}=bG$.
\end{itemize}
\end{corollary}
\begin{proof}
(i)$\Rightarrow$(ii) Since $G'=G$, the relation $\mR$ coincides
with $\sim$. Therefore, it is enough to apply Theorem
\ref{th_taka_1}.

(ii)$\Rightarrow$(iii) This is \cite[Theorem 4.4]{gal_her:i0}.

(iii)$\Rightarrow$(i) It follows from the fact that $G^{xx}$ can
be injected canonically into $G_T$ and the latter can be injected
into $bG$.
\end{proof}
\begin{corollary}\label{co_taka_2}
Let $G$ be a discrete MAP group that is nilpotent of length two,
and such that for each positive integer $n$ there are only
finitely many co-finite normal subgroups $H$ of $G'$ whose index
is less or equal than $n$. Then $G^{xx}\cong G_T$.
\end{corollary}
\begin{proof}
We first verify that, for each $n\in \N$, there is a finite
subgroup $F$ of $U(n)$ such that, given any $D\in \irrng$, we have
$D(G')\subset F$. Indeed, by Lie, Kolchin, Mal'cev's theorem (see
\cite[15.1.1]{robinson}), there is an Abelian subgroup $A$ of
$D(G)$ with finite index, say $m$. That is, $D(g^m)\in A$ for all
$g\in G$. Now, for $x$ and $y$ arbitrarily taken in $G$, since
$G'$ is central in $G$, we have
$(D[x,y])^{m^2}=[D(x),D(y)]^{m^2}=[D(x)^m,D(y)^m]=1$. This yields
$exp~D(G')\leq m^2$ for all $D\in \irrng$. Hence, it suffices to
take $F$ to be the subgroup consisting of all $m^2$-roots of the
identity in $U(n)$.

Thus, fixed any positive integer $n$, we may find $m(n)\in \N$
with $[G':(ker D)\cap G']\leq m(n)$ for all $D\in \irrng$. By
hypothesis, the are only finitely many co-finite normal subgroups
of $G'$ whose index is less or equal than $m(n)$. Therefore, it
follows that $\widehat{G}_n\diagup \mR$ is finite for all $n\in
\N$ and the proof is complete by Theorem \ref{th_taka_1}.
\end{proof}

As a consequence, we obtain the following result due to Poguntke
\cite{poguntke:zwei_klassen,poguntke:chu-dualitat}.
\begin{corollary}\label{co_taka_3}(Poguntke, 1976)
The Heisenberg integral group $H$, satisfies that $H^{xx}\cong
H_T$.
\end{corollary}
\begin{theorem}\label{th_taka_4}
Let $G$ be a discrete MAP group that is an FC group and, for each
positive integer $n$, there are only finitely many co-finite
normal subgroups $H$ of $G'$ such that $G'/H$ accepts faithful
representations into $U(n)$. Then $G^{xx}\cong G_T$.
\end{theorem}
\begin{proof}
Given $n\in \N$ and $D$ arbitrarily chosen in $G^x_n$, we know by
Corollary \ref{co_fc1} that $G'\cap ker D$ is a co-finite normal
subgroup of $G'$. Since, there are a finite number of such
subgroups, we can take $L$ to be the intersection of all them.
Thus, $(G^x_{n})_{|G'}$ can be injected, as a set, into
$(G'/L)^x_n$. Now, the set $\widehat{(G'/L)}_n$ is trivially
finite. This implies that $\widehat{G}_n/\mR$ is finite for all
$n\in \N$. The proof is completed by applying Theorem
\ref{th_taka_1}.
\end{proof}
\begin{corollary}\label{co_taka_5}
Let $G = \sum_{n \in \N} F_n$, where each $F_n\hspace{0.5mm} '$ is
simple and $lim_{n\rightarrow \infty } exp(F_n\hspace{0.5mm}
')=\infty$. Then $G^{xx}\cong G_T$.
\end{corollary}

\section{Some General Remarks}
The results stated in the sections above establish a
classification of countable discrete MAP groups according to their
unitary representation spaces. Let $G$ be a countable discrete MAP
group. From the space of finite dimensional unitary
representations, we have the following cases:

\begin{enumerate}
\item[(i)] There is an integer $n$ and a representation $D\in
G_n^x$ which is faithful. In this case $G$ is a subgroup of
$U(n)$, equipped with the discrete topology.

\item[(ii)] There is an integer $n$ such that $G^x_n$ separates
the points of $G$ but no representation on $G$ is faithful. For
each $D\in G^x_n$, let $N(D)$ be the kernel of $D$. We may form
the group $\prod_{D\in G^x_n} G/N(D)$. Since $G^x_n$ separates
points, the group $G$ may be viewed as a subgroup of $\prod_{D\in
G^x_n} G/N(D)$. Examples of these groups are weak direct sums of
finite groups with bounded exponent equipped with the discrete
topology.

\item[(iii)] There is no integer $n$ such that $G_n^x$ separate
the points of $G$. For example $\sum GF(p_n)$ for various primes
$p_n$ or the Heisenberg integral group. In this case we have that
$G$ is never a Chu group. Indeed, given any integer $n$, there are
infinitely many points which belong to $[G^x_n,V]$ for any $V$
\nbd \ of the identity in $U(n)$. As a consequence, the topology
on $G$ inherited from $G^{xx}$ is not discrete. Hence, the group
is not Chu.
\end{enumerate}
\bigskip

For next result, we need some preparation. Let $G$ and $H$ be
topological groups and let $f$ be a continuous homomorphism from
$G$ to $H$. Then $f$ induces canonically the maps
$f^x_n:H^x_n\longrightarrow G^x_n$, $f^x:H^x\longrightarrow G^x$,
and $f^{xx}:G^{xx}\longrightarrow H^{xx}$. Let $G$ be locally
compact and let $\{F_i\}$ be a collection of compact normal
subgroups such that $F_i\supseteq F_j$ if $i\geq j$, $\cap
F_i=\{1_G\}$, and $G=\underleftarrow{Lim} \ G_i$ (with
$G_i=G/F_i$). Suppose that $H$ has a \nbd \ $V$ of the identity
which contains no small subgroups of $H$ other than $\{1_H\}$.
Then, if $f:G\longrightarrow H$ is any continuous homomorphism,
there is and index $i_0$ such that for $i\geq i_0$ there exists a
continuous homomorphism $f_i:G_i\longrightarrow H$ such that
$f=f_i\circ \phi_i$, where $\phi_i$ is the canonical homomorphism
of $G$ onto $G_i$ (cf. \cite[Chap. 3, p. 294]{bourbaki:GT}).

We shall make use of the following lemma whose proof is left to
the reader.
\begin{lemma}\label{le5.0}
Let $G$ and $H$ be locally compact groups and let
$f:G\longrightarrow H$ a continuous group homomorphism. Then the
diagram
\[
\begin{CD}
G@>f>>   H\\
@V\epsilon_GVV   @VV\epsilon_HV\\     
G^{xx} @>f^{xx}>>  H^{xx}
\end{CD}
\]
\smallskip

\noindent commutes.
\end{lemma}
\begin{proposition}
Let $G=\underleftarrow{Lim} \ G_i$ be a locally compact group
that is the projective limit of 
Chu semi-reflexive groups. Suppose further that each canonical
projection $\phi_i:G\longrightarrow G_i$ is surjective with
compact kernel. Then $G$ is Chu semi-reflexive. If, in addition,
the group $G$ is second countable, then $G$ is a Chu group.
\end{proposition}
\begin{proof}
Since each $G_i$ is Chu and $G$ is embedded into their product, it
follows that $G$ is MAP. We now show that the evaluation map
$\epsilon_{G}$ is surjective. Let $p$ be an element in $G^{xx}$.
For each $i$, $\phi_i^{xx}(p)=\epsilon_{G_i}(x_i)$, with $x_i\in
G_i$. Let $x$ be the element in the product $\prod_{i\in I} G_i$
defined by $(x_i)_{i\in I}$. Using Lemma \ref{le5.0}, it is
readily seen that $x$ belongs to $G$. Now, if $D$ is an arbitrary
element of $G^x$, by the remark above, there is an index $i_0$
such that for each $i\geq i_0$, there exists $D_i\in G_i^x$ such
that $D=D_i\circ \phi_i$. Thus, $p(D)=p(D_i\circ
\phi_i)=p(\phi^x_i(D_i))=\phi_i^{xx}(p)(D_i)=\epsilon_{G_i}(x_i)(D_i)$.
By the way in which $x$ was defined, this yields
$p(D)=\epsilon_{G}(x)(D)$. Therefore $p=\epsilon_{G}(x)$ and
$\epsilon_{G}$ is surjective. Finally, in case that $G$ is second
countable, it suffices to apply Proposition \ref{pr11} to conclude
that $G$ is a Chu group.
\end{proof}

Next example shows that the Chu quasi-dual group $G^{xx}$ need not
be locally compact even for a countable discrete group $G$.
\begin{example}
Let $\{p_n\}$ be an infinite sequence of distinct prime numbers
($p_n>2$), and let $G_n=PSL(2,p_n)$ be the projective special
group of dimension two over the finite filed of order $p_n$. For
each $n$, let $G_{n,m}$ be a copy of $G_n$, for $m=1,2,...$. Let
$G=\sum_{n=1}^{\infty}(\sum_{m=1}^{\infty}G_{n,m})$ with the
discrete topology. The group $G^{xx}$ is not locally compact.
\end{example}
\begin{proof}
Indeed, let $V$ be a small \nbd \ of the identity in $U(n)$. Let
$[G_n^x,V]$ be the \nbd \ in $G^{xx}$ defined by $[G_n^x,V]=\{p\in
G^{xx}: p(G_n^x)\subset V \}$. Observe that
$\epsilon_{G}(\sum_{m=1}^{\infty} G_{l,m})\subset [G_n^x,V]$ if
$l$ is sufficiently large. On the other hand, applying Corollary
\ref{th21}, we know that $\sum_{m=1}^{\infty}G_{l,m}$ is Chu.
Using this and, since the group $G$ projects canonically onto
$\sum_{m=1}^{\infty}G_{l,m}$, it follows that the latter discrete
group is topologically isomorphic to
$\epsilon_{G}(\sum_{m=1}^{\infty}G_{l,m})$. Hence, $G^{xx}$ can
not be locally compact (since no compact set can contain an
infinite discrete subgroup).
\end{proof}

{\bf Acknowledgement}: We would like to thank G\'abor Luk\'acs for
several useful remarks. We also want to thank the referee for
his/her constructive report. They have helped us to improve parts
of this paper.

\newpage

\end{document}